\documentclass{easychair}

\usepackage{amsmath}
\usepackage{tikz}
\usepackage{pgfplots}
\usepackage{booktabs}
\usepackage{hyperref}
\usepackage{xnewcommand}

\newcommand{\obar}[1]{\overline{#1}}
\newcommand{\lit}{\ell}

\newcommand{\approximate}[1]{\hat{#1}}
\newcommand{\approxP}{\approximate{P}}

\newcommand{\approxv}{\approximate{v}}
\newcommand{\approxV}{\approximate{V}}
\newcommand{\approxw}{\approximate{w}}
\newcommand{\approxW}{\approximate{W}}
\newcommand{\approxs}{\approximate{s}}
\newcommand{\round}{\mathit{Rnd}}
\newcommand{\aerror}{\delta}
\newcommand{\digitprecision}{\Delta}

\newcommand{\roundepsilon}{\varepsilon}

\newcommand{\vmin}{v^{-}}
\newcommand{\vmax}{v^{+}}
\newcommand{\interval}[1]{[\![#1]\!]}
\newcommand{\varset}{X}
\newcommand{\dependencyset}{{\cal V}}
\newcommand{\assign}{\alpha}
\newcommand{\modelset}{{\cal M}}

\newcommand{\fexp}{\textsf{exp}}
\newcommand{\ffrac}{\textsf{frac}}

\newtheorem{theorem}{Theorem}
\newtheorem{lemma}{Lemma}

\definecolor{redorange}{rgb}{0.878431, 0.235294, 0.192157}
\definecolor{lightblue}{rgb}{0.552941, 0.72549, 0.792157}
\definecolor{clearyellow}{rgb}{0.964706, 0.745098, 0}
\definecolor{clearorange}{rgb}{0.917647, 0.462745, 0}
\definecolor{mildgray}{rgb}{0.54902, 0.509804, 0.47451}
\definecolor{softblue}{rgb}{0.643137, 0.858824, 0.909804}
\definecolor{bluegray}{rgb}{0.141176, 0.313725, 0.603922}
\definecolor{lightgreen}{rgb}{0.709804, 0.741176, 0}
\definecolor{darkgreen}{rgb}{0.152941, 0.576471, 0.172549}
\definecolor{redpurple}{rgb}{0.835294, 0, 0.196078}
\definecolor{midblue}{rgb}{0, 0.592157, 0.662745}
\definecolor{clearpurple}{rgb}{0.67451, 0.0784314, 0.352941}
\definecolor{browngreen}{rgb}{0.333333, 0.313725, 0.145098}
\definecolor{darkestpurple}{rgb}{0.396078, 0.113725, 0.196078}
\definecolor{greypurple}{rgb}{0.294118, 0.219608, 0.298039}
\definecolor{darkturquoise}{rgb}{0, 0.239216, 0.298039}
\definecolor{darkbrown}{rgb}{0.305882, 0.211765, 0.160784}
\definecolor{midgreen}{rgb}{0.560784, 0.6, 0.243137}
\definecolor{darkred}{rgb}{0.576471, 0.152941, 0.172549}
\definecolor{darkpurple}{rgb}{0.313725, 0.027451, 0.470588}
\definecolor{darkestblue}{rgb}{0, 0.156863, 0.333333}
\definecolor{lightpurple}{rgb}{0.776471, 0.690196, 0.737255}
\definecolor{softgreen}{rgb}{0.733333, 0.772549, 0.572549}
\definecolor{offwhite}{rgb}{0.839216, 0.823529, 0.768627}
\definecolor{medgreen}{rgb}{0.15, 0.6, 0.15}

\definecolor{erd}{RGB}{49,97,160}
\definecolor{mpflow}{RGB}{96,165,200}
\definecolor{mpfmed}{RGB}{161,207,223}
\definecolor{mpfhigh}{RGB}{224,243,248}
\definecolor{mpfilow}{RGB}{215,48,39}
\definecolor{mpfimed}{RGB}{253,174,97}
\definecolor{mpfihigh}{RGB}{254,224,144}
\definecolor{mpq}{RGB}{100,100,100}

\title{Numerical Considerations \\in Weighted Model Counting}

\author{Randal E. Bryant}

\institute{
  Carnegie Mellon University \\
  Pittsburgh, Pennsylvania 15221, USA
}

\authorrunning{R. E. Bryant}
\titlerunning{Numerical Considerations in Weighted Model Counting}

\begin{document}

\maketitle

\begin{abstract}
 Weighted model counting computes the sum of the rational-valued weights
  associated with the satisfying assignments for a Boolean formula,
  where the weight of an assignment is given by the product of the
  weights assigned to the positive and negated variables
  comprising the assignment.  Weighted model counting finds
  applications across a variety of domains including probabilistic reasoning
  and quantitative risk assessment.

  Most weighted model counting programs operate by (explicitly or
  implicitly) converting the input formula into a form that enables
  \emph{arithmetic evaluation}, using multiplication for conjunctions
  and addition for disjunctions.  Performing this evaluation using
  floating-point arithmetic can yield inaccurate results, and it
  cannot quantify the level of precision achieved.  Computing with
  rational arithmetic gives exact results, but it is costly in both
  time and space.

  This paper describes how to combine multiple numeric representations
  to efficiently compute weighted model counts that are guaranteed to
  achieve a user-specified precision.  When all weights are
  nonnegative, we prove that the precision loss of arithmetic
  evaluation using floating-point arithmetic can be tightly bounded.
  We show that supplementing a standard IEEE double-precision
  representation with a separate 64-bit exponent, a format we call
  \emph{extended-range double} (ERD), avoids the underflow and overflow
  issues commonly encountered in weighted model counting.  For problems
  with mixed negative and positive weights, we show that a combination of interval
  floating-point arithmetic and rational arithmetic can achieve the
  twin goals of efficiency and guaranteed precision.  For our
  evaluations, we have devised especially challenging formulas and
  weight assignments, demonstrating the robustness of our approach.

\end{abstract}

\section{Introduction}

Model counting extends traditional Boolean satisfiability (SAT) solving by
asking not just whether a formula can be satisfied, but to compute the
number of satisfying assignments~\cite{gomes:hs:2009}.  Model counting is a challenging
problem---more challenging than the already NP-hard Boolean
satisfiability~\cite{valiant:siam:1979}.

Weighted model counting extends standard model counting by having
rational-valued weights associated with the assignments, and then
computing the sum of the weights of the satisfying assignments.  The
most common variant has weights $w(x)$ and $w(\obar{x})$
assigned to each variable $x$ and its negation $\obar{x}$.  The
weight of an assignment is then the product of the weights for
the positive and negated variables comprising the assignment.
Standard model counting can be seen as a special case of weighted model
counting with all variables and their negations having unit weights: $w(x) = w(\obar{x}) = 1$.

Weighted model counting has applications across a variety of domains,
including probabilistic inference~\cite{chavira:ai:2008,dubray:cp:2024}, Bayesian
inference~\cite{sang:aaai:2005}, probabilistic
planning~\cite{domshlak:jair:2007}, and product line modeling~\cite{sundermann:eme:2023,sundermann:tsem:2024}.  In addition, many of the
applications of decision diagrams (DDs) for
combinatorics~\cite{knuth:bdd:2011}, quantitative risk
assessment~\cite{andrews:ieeetr:2000,groen:ress:2006,hardy:ieeer:2007,xing:wiley:2015,xing:amm:2025},
Bayesian inference~\cite{minato:ijcai:2007},
optimization~\cite{bergman:book:2016}, and
product line modeling~\cite{andersen:jair:2010,benavides:is:2010} are, at their core, applications of
weighted model counting for discrete functions represented as decision diagrams.

Despite the intractability, a variety of weighted model counting
programs have been developed that work well in practice.  They
generally fall into two categories~\cite{shaw:kr:2024}. \emph{Top-down} programs
recursively branch on the variables of a formula, performing unit
propagation and conflict analysis similar to CDCL SAT solvers.  Most
of these programs operate as \emph{knowledge compilers}, converting
the input Boolean formula into a restricted form that enables efficient
weighted and unweighted counting~\cite{darwiche:aaai:2002,darwiche:ecai:2004,lagniez:ijcai:2017,muise:cai:2012,oztok:cp:2014,sharma:ijcai:2019}.
Others apply \emph{bottom-up} approaches, including ones using
multi-terminal BDDs~\cite{dudek:aaai:2020,dudek:sat:2021}.  In both
cases, the strategy is to convert the formula into a form for which
weighted model counting becomes tractable.

Weighted model counting can be computationally intensive.  In
experimental results described in this paper, some evaluations 
require over one billion arithmetic operations.
Floating-point arithmetic can provide the needed level of performance,
but the computed values are often either too small or too
large in magnitude to encode with standard floating-point
representations.  In addition, the rounding errors
introduced by floating-point computations can lead to results
that bear little relation to the actual values.  In general, even when the results are accurate,
floating-point evaluation cannot quantify the level of precision achieved.

Absolute precision can be guaranteed by performing the
computations with a rational-arithmetic software package~\cite{knuth:rational:1981}, such as the
MPQ library within the GNU Multiprecision Arithmetic
Library (GMP)~\cite{granlund:gmp:2015}.  It represents a rational number $v$ as a
pair of multiprecision integers $p$ and $q$ with $v = p/q$.
MPQ
can compute the exact rational values of all multiplication and
addition operations, yielding an exact weighted count.  Unfortunately,
both the space and the time required for storing and manipulating these numbers can be very
large.
In this paper, for example, we report experiments requiring over one gigabyte
to store the arguments and result of a single addition operation.
For most applications, rational arithmetic provides more precision than is required.
It would be preferable to have a floating-point
representation, with the ability to set and achieve a level of precision suitable for a given application.

This paper describes how to combine multiple representations to
compute weighted model counts that are guaranteed to achieve a
user-specified precision, enabling a tradeoff between precision and
computation time.  First, we consider the case where the weights
$w(x)$ and $w(\obar{x})$ for all variables $x$ are nonnegative, and
where the values are computed over a decision-DNNF Boolean
formula~\cite{beame:uai:2013,huang:jair:2007}.  We prove under these
restrictions that the degradation of precision caused by rounding
errors will be bounded by the logarithm of the number of variables in
the formula.  In practical terms, this implies that floating-point
arithmetic can be fully trusted in these cases.

Our experiments show that the standard IEEE double floating-point
representation is prone to underflow and overflow when performing
weighted-model counting.  We have developed the \emph{Extended-Range
Double} (ERD) floating-point library to overcome this limitation by
augmenting a standard IEEE double with a separate exponent field
stored as a 64-bit signed number.  To achieve higher precision, we use
the MPF software floating-point library within GMP with fraction sizes
$p \in \{64, 128, 256\}$, depending on the target precision.
The 64-bit exponent fields of ERD and MPF suffice to represent
the full range of values in weighted model counting.

This result has broad applicability.
For many applications
of weighted model counting, the weights are probabilities between
$0.0$ and $1.0$, or they are unit weights.  For these applications,  negative weights are never encountered.
Furthermore, most top-down and bottom-up weighted model
counters either explicitly or implicitly operate on decision-DNNF
representations~\cite{beame:uai:2013}.  
Decision diagrams with binary branching structure
also have direct translations into decision-DNNF formulas~\cite{huang:jair:2007,oztok:cp:2014}.

To extend this capability to less restricted classes of Boolean
formulas and to decision diagrams with nonbinary branching
structures~\cite{darwiche:ijcai:2011,srinivasan:iccad:1990}, we
present a method for computing an integer-valued error bound prior to arithmetic
evaluation.  This bound can guide the selection of the floating-point fraction size to
achieve the desired precision.

For formulas with \emph{mixed} negative and positive weights,
we show experimentally that floating-point
arithmetic suffices for the common ways weight assignments are generated in benchmark evaluations.
On the other hand, we
describe a strategy for generating random weight assignments that
often causes floating-point arithmetic to yield erroneous results.
We also
demonstrate a family of
formulas where no bounded-precision numerical representation will
suffice.  Rational arithmetic provides the only option in such cases.

We address the lack of certainty in floating-point
evaluation by introducing \emph{interval} floating-point
arithmetic~\cite{hickey:jacm:2001} using the MPFI software
library~\cite{revol:rc:2005}. With this library, values are
represented by closed intervals, written $\interval{\vmin, \vmax}$, such that the
true value $v$ satisfies $\vmin \leq v \leq \vmax$, and both $\vmin$ and $\vmax$ are
represented in floating point.  The result of every operation is an
interval that is guaranteed to include the true result value, as long
as the argument intervals include their true values~\cite{hickey:jacm:2001,muller:hfpa:2018}.
When an interval must be converted to a single value, the floating-point number nearest the midpoint $(\vmin+\vmax)/2$ is chosen.
We show
experimentally that the intervals maintained during the computations
of weighted model counting are generally tight enough to provide
useful precision guarantees.

Putting these together, we present experimental results for a program
that employs a hybrid strategy to compute the weighted count of a decision-DNNF formula
generated by the D4 knowledge compiler~\cite{lagniez:ijcai:2017}.
The user specifies a target precision $D$, measured in decimal digits, as defined in Section~\ref{sect:background:numbers}.
When all weights are nonnegative, it uses either our ERD representation or MPF with an appropriate fraction size to
perform floating-point computations, relying on our precision
guarantee.  For mixed weights, it performs multiple levels
of interval computation with MPFI, using increasing precision.  If these evaluations fail to guarantee the target precision,
it resorts to rational arithmetic using MPQ\@.  The overall effect
is to achieve the twin goals of efficiency and guaranteed precision.
Although we only present experimental results for D4, similar results will hold for other top-down and bottom-up weighted model counters, as well as
for numerical computations on decision diagrams.

This work is motivated by both application need and technical
opportunity.  On the need side, there is evidence that the standard of
precision for current weighted model counters is low.  In the 2020
weighted model counting competition, a count was considered correct if
it was within $10\%$ of a precomputed result~\cite{fichte:jea:2020},
corresponding to decimal precision $D=1$.  That threshold has been
tightened to $0.1\%$ (decimal precision $D=3$) in more recent
years~\cite{hecher:mc:2024}.  Such low precision may suffice for some
applications, but it is significantly below the level achieved by
other numerical programs.  On the opportunity side, our work
demonstrates the ability to achieve target precisions ranging from
$D=1$ to $D=70$, while generally avoiding the high cost of rational
arithmetic.

We see this work as going beyond satisfying the needs of current
applications of weighted model counting to create a robust approach
that will handle future applications.  To test robustness, we have
devised formulas and weight assignments that present challenging cases
for numerical accuracy.  We show even these cases can be handled by an
appropriate combination of numerical representations.

Regarding previous work, most recent work on estimating the error
caused by floating-point rounding \emph{a priori} focuses on getting precise bounds
and supporting a variety of
operations, but with less concern about scalability~\cite{becker:fmcad:2016,magron:toms:2017,solovyev:toplas:2018}.
By contrast, we only seek
loose bounds and only when multiplying and
adding nonnegative numbers.  On the other hand, we must be able to scale to evaluations consisting of billions of operations.
Consequently, we reach back to
more historic work~\cite{wilkinson:nm:1960,wilkinson:rounding:1964}.
We have not seen any investigation of the numerical properties of weighted model counting, and especially
the tight error bounds that can be obtained for the arithmetic evaluation of decision-DNNF formulas.
We also have not seen any systematic studies on the performance of interval or rational arithmetic for
weighted model counting.

Sections~\ref{sect:background:boolean} and
\ref{sect:background:numbers} of this paper cover background material
in Boolean formulas, weighted model counting, numerical error, and
numeric representations.  Section~\ref{sect:nonneg} covers the case of
nonnegative weights, with our main theoretical result, a means of
computing error bounds for more general formulas and decision
diagrams, and an experimental validation.  Section~\ref{sect:neg}
describes the challenges that negative weights can present, but also
experimental results showing that floating-point arithmetic suffices
in many cases.  It describes the use of interval arithmetic of
increasing precision, along with rational arithmetic, to reliably
handle challenging benchmarks.

Section~\ref{sect:hybrid} describes and evaluates our hybrid approach.
Section~\ref{sect:erd} describes our method for extending the range of the IEEE double representation.
Section~\ref{sect:conclusion}  presents some concluding remarks.

\section{Boolean Formulas and Weighted Model Counting}
\label{sect:background:boolean}

We consider Boolean formulas over a set of variables $\varset$ in
\emph{negation normal form}, where negation can only be applied to the
variables.  We refer to a variable $x$ or its negation $\obar{x}$ as a
\emph{literal}.  We use the symbol $\lit$ to indicate an arbitrary
literal.  The set of all formulas is defined recursively to consist of
literals, \emph{conjunctions} of the form $\phi_1 \land \phi_2$, and
\emph{disjunctions} of the form $\phi_1 \lor \phi_2$.  The set of
variables occurring in formula $\phi$ is denoted
$\dependencyset(\phi)$.  Typically, a formula is represented as a
directed acyclic graph, allowing a sharing of subformulas.  We
therefore define the \emph{size} of a formula to be the number of
unique subformulas.

A (total) assignment is a mapping $\assign \colon \varset \rightarrow
\{0, 1\}$.  Assignment $\assign$ is said to be a \emph{model} of
formula $\phi$ if the formula evaluates to $1$ under that assignment.
The set of models of a formula $\phi$ is written $\modelset(\phi)$.
We can also consider an assignment to be a set of literals, where $x \in \assign$
when $\assign(x) = 1$, and $\obar{x} \in \assign$ when
$\assign(x) = 0$, for each variable $x \in \varset$.

\emph{Weighted model counting} is defined in terms of a \emph{weight
  assignment} $w$, associating rational values $w(x)$ and
$w(\obar{x})$ with each variable $x \in \varset$.
The weight of an
assignment $\assign$ is then defined to be the product of its literal weights, and the weight
of a formula is the sum of the weights of its satisfying assignments:
\begin{eqnarray}
  w(\phi) & = & \sum_{\assign \in \modelset(\phi)} \;\;\prod_{\lit \in \assign} w(\lit) \label{eqn:wmc}
\end{eqnarray}

Computing the weighted count of an arbitrary formula is thought to be intractable.  However, it becomes
feasible when the formula is in \emph{deterministic decomposable} negation-normal form (d-DNNF):
\begin{enumerate}
\item The formula is in negation-normal form.  
\item All conjunctions are \emph{decomposable}~\cite{darwiche:jacm:2001,darwiche:jair:2002}.  That is, every subformula $\phi'$ of the form $\phi' = \phi_1 \land \phi_2$
  satisfies $\dependencyset(\phi_1) \cap \dependencyset(\phi_2) = \emptyset$.
\item All disjunctions are \emph{deterministic}~\cite{darwiche:jair:2002,darwiche:jancl:2001}.  That is, every subformula $\phi'$ of the form $\phi' =\phi_1 \lor \phi_2$ satisfies
  $\modelset(\phi_1) \cap \modelset(\phi_2) = \emptyset$.
\end{enumerate}
As an important subclass of d-DNNF, a formula is said to be in 
\emph{decision decomposable} negation-normal form (decision-DNNF)~\cite{huang:jair:2007} when every occurrence of a disjunction has the form 
$(x \land \phi_1) \lor (\obar{x} \land \phi_2)$ for some variable $x$, referred to as the \emph{decision variable}.

There are several ways to compute the weighted count of d-DNNF formula $\phi$, all
based on an \emph{arithmetic evaluation} of $\phi$ to compute a value $W(\phi)$:
\begin{enumerate}
\item Each literal $\lit$ is a assigned a rational value $W(\lit)$.
\item Each subformula $\phi' = \phi_1 \land \phi_2$ is evaluated as $W(\phi') = W(\phi_1) \cdot W(\phi_2)$.
\item Each subformula $\phi' = \phi_1 \lor \phi_2$ is evaluated as $W(\phi') = W(\phi_1) + W(\phi_2)$.
\end{enumerate}
The number of arithmetic operations in this evaluation is linear in the size of the formula.

The following methods use arithmetic evaluation to compute a weighted model count $w(\phi)$ of formula $\phi$ for weight assignment $w$:
\begin{enumerate}
\item If the weight assignment satisfies $w(x) + w(\obar{x}) = 1$  for every variable $x$,
  then by letting $W(\lit) = w(\lit)$ for each literal $\lit$, the arithmetic evaluation $W(\phi)$ will yield the weighted model count $w(\phi)$.
\item Formula $\phi$ is said to be \emph{smooth} if every disjunction $\phi_1 \lor \phi_2$ satisfies
  $\dependencyset(\phi_1) = \dependencyset(\phi_2)$~\cite{darwiche:jair:2002,darwiche:jancl:2001}.  For a smooth formula, 
by letting $W(\lit) = w(\lit)$ for each literal $\lit$, the arithmetic evaluation $W(\phi)$ will yield the weighted model count $w(\phi)$.
\item If the weight assignment satisfies $w(x) + w(\obar{x}) \not = 0$ for every variable $x$,
  we can apply \emph{rescaling}, first computing $s(x) = w(x) + w(\obar{x})$ for each variable $x$
  and letting $W(x) = w(x)/s(x)$ and $W(\obar{x}) = w(\obar{x})/s(x)$.  
  Following the arithmetic evaluation, the weighted count is computed as:
  \begin{eqnarray}
w(\phi) &=& W(\phi)\; \cdot \;  \prod_{x\in\varset} s(x)  \label{eqn:rescale}
  \end{eqnarray}
\end{enumerate}

An arbitrary formula can be smoothed by inserting \emph{smoothing terms} of the form $x \lor \obar{x}$~\cite{darwiche:jancl:2001}.
For example,
  a disjunction $\phi_1 \lor \phi_2$ having $x \in \dependencyset(\phi_1)$ but
  $x \not\in \dependencyset(\phi_2)$ is rewritten as $\phi_1 \lor [(x \lor \obar{x}) \land \phi_2]$.
  Adding smoothing terms can expand the size of a formula significantly, and it can be time consuming.
  More precisely, for $n = |X|$, and a formula with $m$ unique subformulas, it can require time $\Theta(m\cdot n)$ and increase the formula size by a factor of $n$.
  Some restricted formula classes allow more space- and time-efficient smoothing~\cite{shih:neurips:2019}, including those arising from decision
  diagrams with totally ordered variables~\cite{bryant:ieeetc:1986,minato:sttt:2001}.
  However, the required properties do not hold for the formulas generated by most weighted model counters.
  
  These three methods can be combined by rescaling some variables,
  inserting smoothing terms for others, and taking no action for the rest.
  For example,
  for any variable $x$ having $w(x) + w(\obar{x}) = 0$, we can insert
  smoothing terms, while applying rescaling for other variables $y$ such that $w(y) + w(\obar{y}) \not= 1$.

\section{Approximations and Numeric Representations}
\label{sect:background:numbers}

When approximating rational number $v$ with value $\approxv$, we
define the \emph{approximation error} $\aerror[\approxv, v]$
as the relative error when $v \not = 0$ and as requiring an exact representation when $v = 0$:
\begin{eqnarray}
\aerror[\approxv, v] & = & \left\{ \begin{array}{ll}
  \frac{|\approxv - v|}{|v|}  & v \not = 0\\
  0 & v  = \approxv = 0\\
  1 & v = 0 \; \textrm{and} \; \approxv \not = 0
  \end{array} \right. \label{eqn:approx:error}
\end{eqnarray}
This value will equal 0 when $\approxv=v$, and it will be greater for weaker approximations.
Observe that approximating a nonzero value with zero yields a high error:  $\aerror[0, v] = 1$ when $v\not=0$.

The \emph{decimal precision} expresses the quality of an approximation by the number of significant digits in its decimal representation:
\begin{eqnarray}
\digitprecision(\approxv, v) & = & \max(0, -\log_{10} \aerror[\approxv, v]) \label{eqn:digitprecision} 
\end{eqnarray}
This value will range from $0$ for a poor approximation, up to $+\infty$ when $\approxv=v$.

We consider floating-point numbers of the form
\begin{eqnarray}
v & = & (-1)^s \; \times \; f \; \times 2^{e} \label{eqn:floating-point}
\end{eqnarray}
where:
\begin{itemize}
\item Sign bit $s$ equals $0$ for nonnegative numbers and $1$ for negative numbers
\item Fraction $f$ is encoded as a $p$-bit binary number with an implicit binary point on the left.  That is $0 \leq f \leq 1-2^{-p}$.
\item Exponent $e$ is an integer, possibly with some limitation on its range.
\end{itemize}
As examples, consider two different floating-point formats:
\begin{itemize}
\item The IEEE~754 Double format uses a slightly different
  representation, but it maps to the notation of
  Equation~\ref{eqn:floating-point} with $p=53$ and an exponent range
  of $-1021 \leq e \leq 1024$~\cite{overton:siam:2001}.
  Unfortunately, the small exponent range (giving a magnitude range
  for nonzero numbers of around $10^{\pm 308}$) limits the suitability
  of this representation for weighted model counting.  For example, as
  part of the evaluation of weighted model counting when all weights
  are nonnegative, described in Section~\ref{sect:nonneg}, we computed
  the counts for 1000 combinations of formula and randomly-generated
  weight assignment using double-precision arithmetic.  Fully 628 of
  the evaluations failed due to values exceeding the exponent range,
  with 419 overflowing to infinity and 209 underflowing to zero.  For
  the original weight assignments provided with the 100 formulas
  evaluated, 45 of them failed with double-precision arithmetic, with
  5 overflowing and 40 underflowing.  To counter this deficiency, we
  have implemented a floating-point library using an
  \emph{Extended-Range Double} (ERD) numerical representation,
  augmenting an IEEE Double with a 64-bit signed exponent, as
  discussed in Section~\ref{sect:erd}.

\item The MPF software floating-point library allows the value of $p$
  to be set to any multiple of 64.  We use configurations with $p$ equal to $64$, $128$, and $256$,
  referring to these as ``MPF-64,'' ``MPF-128,'' and ``MPF-256.'' On
  most 64-bit architectures, MPF represents the exponent as a 64-bit
  signed number.  This provides an ample exponent range, giving a
  magnitude range of over $10^{\pm 10^{18}}$.
  For example, the weighted model count for a tautology with $n$ variables,
  where all literals are assigned weight $w$, equals
  $2^n\cdot w^n = (2w)^n$.  Consider literal weight $w=10^{1000}$,
 far larger than can even be represented as an IEEE-754 Double, and let $n$ equal one trillion, over five orders of magnitude
 larger than the largest formulas solved by current weighted model counters.  The weighted count is
 $(2 \times 10^{1000})^{10^{12}} \approx 10^{10^{15.00013}}$.  In the other direction, the conjunction of one trillion variables, each having a weight
 of $10^{-1000}$ has a weighted count of $10^{-10^{15}}$.  Even these extreme values are
 well within the range of the MPF representation.
\end{itemize}

From this we can conclude: 1) IEEE Double can be used when a
fraction size of $p=53$ suffices, and its range limitation can be overcome
through our ERD representation,
and 2) MPF can use fraction sizes that provide very high precision.  We
assume for the remainder of this paper that all floating-point
computations can be performed without underflow or overflow.

\begin{table}
  \caption{Bounds on Round-Off Errors for Different Floating Point Representations.  The lower bound on
  weighted model counting (rightmost column) holds for formulas with $n \leq 10^{7}$ variables when all weights are nonnegative.}
  \label{tab:precision}
  \begin{center}
  \begin{tabular}{lrrrr}
    \toprule
    \multicolumn{1}{c}{Bound Type} & & \multicolumn{1}{c}{Upper} & \multicolumn{1}{c}{Lower} & \multicolumn{1}{c}{Lower} \\
    \multicolumn{1}{c}{} & \multicolumn{1}{c}{$p$} & \multicolumn{1}{c}{$\roundepsilon$} & \multicolumn{1}{c}{$\digitprecision(\round(v), v)$} & \multicolumn{1}{c}{$\digitprecision(\approxw(\phi), w(\phi))$} \\
    \midrule
    IEEE Double / ERD &  53 & $1.11 \times 10^{-16}$ & $15.95$    & $8.11$ \\
    MPF-64      &  64 & $5.42 \times 10^{-20}$   & $19.27$  & $11.42$ \\
    MPF-128      &  128 & $2.94 \times 10^{-39}$ &  $38.53$ & $30.69$  \\
    MPF-256      &  256 & $8.64 \times 10^{-78}$ & $77.06$ & $69.22$ \\
    \bottomrule
  \end{tabular}
  \end{center}
\end{table}

When encoding rational number $v$ in floating point, its value must be
rounded to a value $\round(v)$.  Doing so can introduce rounding
error~\cite{knuth:fp:1981,muller:hfpa:2018}.  Letting $\roundepsilon =
2^{-p}$, we can assume that $\aerror[\round(v), v] \leq
\roundepsilon$, and that $\digitprecision(\round(v), v) \geq p \,
\log_{10} 2$.
The third column of
Table~\ref{tab:precision} lists the bounds on $\roundepsilon$ for the
four different floating-point representations considered, while the fourth column lists the bounds on $\digitprecision$.

Floating-point arithmetic is implemented in such a way that any
operation effectively computes an exact result and then rounds it to
encode the result as a floating-point value.  The maximum error from a sequence of operations
therefore tends to accumulate in multiples of $\roundepsilon$.
This yields error bounds of the form $\aerror[\approxv, v] \leq t\,\roundepsilon$,
which we  refer to as having at most $t$ units of
rounding error.

For interval $\interval{\vmin, \vmax}$, we define the interval approximation error $\aerror(\interval{\vmin, \vmax})$ as
\begin{eqnarray}
\aerror(\interval{\vmin, \vmax}) & = & \left\{ \begin{array}{ll}
  \frac{\vmax - \vmin}{\min(|\vmin|, |\vmax|)}  & 0 \not \in \interval{\vmin, \vmax}\\[0.8em]
  0 & \vmin = \vmax = 0 \\
  1 & 0 \in \interval{\vmin, \vmax} \;\; \textrm{and} \;\; \vmin < \vmax
  \end{array} \right. \label{eqn:interval:error}
\end{eqnarray}
For any values $\approxv, v \in \interval{\vmin, \vmax}$, we can see that
$\aerror[\approxv, v] \leq \aerror(\interval{\vmin, \vmax})$.
We then define the decimal precision of the interval as:
\begin{eqnarray}
\digitprecision(\interval{\vmin, \vmax}) & = & \max[0, -\log_{10} \aerror(\interval{\vmin, \vmax})] \label{eqn:interval:digitprecision} 
\end{eqnarray}
This value will range from $0.0$ for a very large interval, relative to the magnitudes of its endpoints, to $+\infty$ when the interval is tight with $\vmin = \vmax$.

\section{Only Nonnegative Weights}
\label{sect:nonneg}

Here we evaluate how rounding errors accumulate via a series of
arithmetic operations when all arguments are nonnegative.
That is, assume the exact arguments $v$ and $w$ for each operation satisfy $v \geq 0$ and $w \geq 0$.
Rounding never causes a nonnegative number to become negative, and therefore 
the approximations $\approxv$ of  $v$ and $\approxw$ of $w$ satisfy $\approxv \geq 0$ and $\approxw \geq 0$.
None of the operations multiplication, addition, or division yield negative results when their arguments are nonnegative.
We can therefore
assume that all actual and approximate values under consideration are
nonnegative.

Our analysis builds on historic work for bounding the error produced
by a series of floating-point
multiplications~\cite{muller:hfpa:2018,rump:bit:2015,wilkinson:nm:1960,wilkinson:rounding:1964} or
additions~\cite{higham:siam:1993}.  Our formulation considers combinations of multiplication and addition, and it weakens the error bound to simplify the analysis.
It
applies only when all arguments are nonnegative.

Suppose for nonnegative values of $v$ and $w$ and nonnegative values $s$ and $t$, we have
$\aerror[\approxv, v] \leq s\, \roundepsilon$ and
$\aerror[\approxw, w] \leq t\, \roundepsilon$, respectively.  Assume also that we have $v = 0$ if and only if $\approxv = 0$, and similarly from $w$ and $\approxw$.
The bounds can be expanded according to (\ref{eqn:approx:error}) as
$(1-s\,\roundepsilon)\, v \leq \approxv \leq (1+s\,\roundepsilon)\, v$ and
$(1-t\,\roundepsilon)\, w \leq \approxw \leq (1+t\,\roundepsilon)\, w$.

\subsection{Multiplication}

Assume that $v > 0$ and $w > 0$ and consider the effect of multiplying their approximations
$\approxv$ and $\approxw$.  To simplify the analysis, let us impose as an additional constraint that $s\,t \leq 1/\roundepsilon$.
The product $\approxv \cdot \approxw$  satisfies
$\approxv \cdot \approxw \leq (v\cdot w) [1 + (s+t)\,\roundepsilon + s\,t\,\roundepsilon^2]$, and we can use the additional constraint to replace $s\,t\,\roundepsilon^2$ by $\roundepsilon$,
giving
$\approxv \cdot \approxw \leq (v\cdot w) [1 + (s+t+1)\,\roundepsilon]$.
In the other direction, $\approxv \cdot \approxw$ satisfies
$(v\cdot w) [1 - (s+t)\,\roundepsilon + s\,t\,\roundepsilon^2] \leq \approxv \cdot \approxw$.  We can drop the term
$s\,t\,\roundepsilon^2$ to give
$(v\cdot w) [1 - (s+t)\,\roundepsilon] \leq \approxv \cdot \approxw$.  These two bounds guarantee that
$\aerror[\approxv \cdot \approxw, v \cdot w] \leq (s+t+1)\,\roundepsilon$.
Rounding this result can introduce an additional error of at most $\roundepsilon$, and therefore
$\aerror[\round(\approxv \cdot \approxw), v \cdot w] \leq (s+t+2)\,\roundepsilon$.

When $v=0$ (respectively, $w=0$), we will have $\approxv=0$ (resp., $\approxw = 0$) and therefore $v \cdot w = \approxv \cdot \approxw = 0$.
We can therefore state that  for any nonnegative values of $v$ and $w$, the three conditions $\aerror[\approxv, v] \leq s\,\roundepsilon$,
$\aerror[\approxw, w] \leq t\,\roundepsilon$, and
$s\,t \leq 1/\roundepsilon$, imply that
$\aerror[\round(\approxv \cdot \approxw), v \cdot w] \leq (s+t+2)\,\roundepsilon$.

Thus, for values of $s$, $t$, and $\roundepsilon$ satisfying our
additional constraint, a multiplication operation, at most, propagates
the sum of the errors of its arguments, and it adds two units of
rounding error.

\subsection{Addition}

When positive values $v$ and $w$ are added, their approximations  $\approxv$ and $\approxw$ satisfy
$(v + w) (1 - r\,\roundepsilon) \leq \approxv + \approxw \leq (v + w) (1 + r\,\roundepsilon)$, where
$r = (s\,v + t\,w)/(v+w)$.
That is, the resulting error bound $r$ is a weighted average
of those of its arguments.  For all values of $v$ and $w$, $r$ cannot exceed the maximum of $s$ and $t$.
Rounding the sum can add at most one unit of rounding error, and so we have
$\aerror[\round(\approxv + \approxw), v + w] \leq (\max(s,t)+1)\,\roundepsilon$.

If $v = 0$ (respectively, $w = 0$), we will have $\approxv = 0$ (resp., $\approxw = 0$), and therefore $\round(\approxv + \approxw) = \approxw$ (resp., $= \approxv$).
We can therefore state that for any nonnegative values of $v$ and $w$, the two conditions $\aerror[\approxv, v] \leq s\,\roundepsilon < 1$ and
$\aerror[\approxw, w] \leq t\,\roundepsilon < 1$ imply that 
$\aerror[\round(\approxv + \approxw), v + w] \leq (\max(s,t)+1)\,\roundepsilon$.

Thus, an addition operation, at most, propagates the maximum error of its arguments, and it adds one unit of rounding error.

\subsection{Evaluating a Decision-DNNF Formula}
\label{sect:error:formula}

Suppose we use floating-point arithmetic to compute the sums and
products in an arithmetic evaluation of a decision-DNNF formula $\phi$.
We assume that the value $W(\lit)$ for each literal $\lit$ is
represented by a floating-point number $\approxW(\lit)$ such that
$\aerror[\approxW(\lit), W(\lit)] \leq \roundepsilon$.  In practice,
this implies that rescaling must use rational arithmetic to
compute exact representations of $s(x)$, $w(x)/s(x)$, and
$w(\obar{x})/s(x)$ for each variable $x$, so that only one unit of
rounding error is introduced when representing each value $W(\lit)$.
We can then bound the error of the computed value $W(\phi)$ as follows:
\begin{lemma}
  The arithmetic evaluation of a decision-DNNF formula $\phi$  having $|\dependencyset(\phi)| = n$, with
  $n \leq 1/(2\sqrt{\roundepsilon})$
  using floating-point arithmetic,
  and where all literals $\ell$ satisfy $W(\lit) \geq 0$,
  will yield an approximation $\approxW(\phi)$ satisfying
  $\aerror[\approxW(\phi), W(\phi)] \leq (4n-2)\,\roundepsilon$.
  \label{lemma:approx:pos}
\end{lemma}

The proof of this lemma proceeds by induction on the structure of $\phi$:
\begin{enumerate}
\item For literal $\lit$ with weight $W(\lit)$, its approximation $\approxW(\lit)$ satisfies
$\aerror[\approxW(\lit), W(\lit)] \leq \roundepsilon$, which is within the error bound of $(4n-2)\,\roundepsilon$ for $n=1$.
\item For conjunction $\phi$ of the form $\phi_1 \land \phi_2$, there must be some $k$, with $1 \leq k < n$, such that $|\dependencyset(\phi_1)| = k$
  and $|\dependencyset(\phi_2)| = n-k$.
\begin{enumerate}
\item Let us first test whether the requirement that $n \leq 1/(2\sqrt{\roundepsilon})$ 
  guarantees that the conditions on $s$, $t$, and $\roundepsilon$ required for the multiplication bound  hold.
For $s = 4k-2$ and $t = 4(n-k)-2$
    we require that $s\,t \leq 1/\roundepsilon$.  We can see that $s\,t \leq 16\,k\,(n-k)$. This quantity will be maximized when $k = n/2$,
    and therefore $s \, t \leq 4\,n^2$.
    Given our limit on $n$ with respect to $\roundepsilon$, we have $s \, t \leq 1/\roundepsilon$.
  \item
    We can also see that if $n \leq 1/(2\sqrt{\roundepsilon})$, then both
    $k \leq 1/(2\sqrt{\roundepsilon})$ and  $n-k \leq 1/(2\sqrt{\roundepsilon})$.
  \item
We can therefore assume by induction that 
 $\aerror[\approxW(\phi_1), W(\phi_1)] \leq (4 k-2) \,\roundepsilon$
  and also that $\aerror[\approxW(\phi_2), W(\phi_2)] \leq (4 (n-k)-2) \,\roundepsilon$.  Their product, after rounding
  will satisfy 
$\aerror[\approxW(\phi_1 \land \phi_2), W(\phi_1 \land \phi_2)] \leq [(4 k -2) + (4 (n-k) -2) + 2]\,\roundepsilon = (4n-2) \,\roundepsilon$.
\end{enumerate}
\item For disjunction $\phi$ of the form
  $\phi = (x \land \phi_1) \lor (\obar{x} \land \phi_2)$, let us use the notation $\lit_1 = x$ and $\lit_2 = \obar{x}$
  and consider the two subformulas $\lit_i \land \phi_i$ for $i \in \{1,2\}$.
  Since all products are decomposable, we must have $x \not \in \dependencyset(\phi_i)$,
  and therefore $|\dependencyset(\phi_i)| \leq n-1$.  We can also see that the condition 
  $n \leq 1/(2\sqrt{\roundepsilon})$ implies that   $n-1 \leq 1/(2\sqrt{\roundepsilon})$.
  By induction, we can therefore assume that
  $\aerror[\approxW(\phi_i), W(\phi_i)] \leq (4(n-1)-2) \,\roundepsilon = (4n-6)\,\roundepsilon$.  Rounding the literal weights will yield
  $\aerror[\approxW(\lit_i), W(\lit_i)] \leq \roundepsilon$.  Let $v_i$ denote the product $W(\lit_i) \cdot W(\phi_i)$ for $i \in \{1,2\}$.
  Its rounded value will satisfy
  $\aerror[\approxv_i, v_i] \leq (4n-3) \,\roundepsilon$.  Summing $\approxv_1$ and $\approxv_2$ and rounding the result will therefore give
  an approximation $\approxW(\phi)$ to $W(\phi) = v_1 + v_2$ with
$\aerror[\approxW(\phi), W(\phi)] \leq (4n-2)\,\roundepsilon$.  
\end{enumerate}

Observe that this proof relies on the decomposability of the
conjunctions to bound the error induced by multiplication operations.
It relies on the decision structure of the formula only to bound the
depth of the additions.  It does not rely on the formula being deterministic.

In the event of rescaling, we must also consider the error introduced
when computing the product $P = \prod_{x\in\varset} s(x)$.  We assume that
each term $s(x)$ is represented by a floating-point value
$\approxs(x)$ such that $\aerror[\approxs(x), s(x)] \leq
\roundepsilon$.  In practice, this requires using rational arithmetic
to represent $w(x)$ and $w(\obar{x})$ and to compute their sum.  The
only error introduced will then be when converting the sum into
a floating-point representation.

We can then bound the error of the product as
\begin{lemma}
  The computation of the product $P = \prod_{x\in\varset} s(x)$ having
$|\varset| = n$, with $n \leq 1/(2\sqrt{\roundepsilon})$ using floating-point arithmetic
and where all literals $\ell$ satisfy $s(\lit) \geq 0$,
will yield an approximation $\approxP$ satisfying
  $\aerror[\approxP, P] \leq (3n-2)\,\roundepsilon$.
  \label{lemma:approx:product}
\end{lemma}

The proof of this lemma proceeds much like that for Lemma~\ref{lemma:approx:pos}.  We assume an arbitrary association of the subproducts, and so $P$ can be computed as
$P_1 \cdot P_2$, where $P_1$ is the product of $k$ elements and $P_2$ is the product of $n-k$ elements, with $1 \leq k < n$.
The smaller coefficient of $3$ arises due to the lack of addition operations.

Combining the two lemmas, we can state the following result about weighted model counting when all weights are nonnegative:
\begin{theorem}
  \label{thm:approx:pos}
Computing the weighted model count of
a decision-DNNF formula $\phi$, where all literal weights are nonnegative, with $|\dependencyset(\phi)| = n$, and using floating-point arithmetic with a $p$-bit fraction, such that $\log_2 n \leq p/2-1$
will yield an approximation $\approxw(\phi)$ to the true weighted count $w(\phi)$, such that
\begin{eqnarray}
\digitprecision(\approxw(\phi), w(\phi)) & \geq & p \cdot \log_{10}2 - \log_{10}n - c\label{eqn:precision:wmc}
\end{eqnarray}
where $c = \log_{10} 7$ when rescaling is required and $c = \log_{10} 4$ when no rescaling is required.
\end{theorem}

Let us examine the practical implications of this theorem.  Assume we are given a decision-DNNF formula over $n$
variables and wish to compute its weighted model count via rescaling, where all weights are nonnegative, with a decimal
precision of at least $D$.  We can do so using a floating-point precision $p$ that satisfies the following two conditions:
\begin{eqnarray}
  p & \geq & 2(1 + \log_2 n) \label{eqn:pmin:linear} \\
  p & \geq & D \cdot \log_2 10 +\log_2 n + 2.9 \label{eqn:pmin:roundoff}
\end{eqnarray}
For example, no formula from the 2024 weighted model counting
competition had more than 10 million variables, and so we can assume
$\log_2 n \leq 23.3$.  Equation~\ref{eqn:pmin:linear} then requires $p \geq 48.6$.
Using $p = 49$, Equation~\ref{eqn:precision:wmc} then guarantees digit precision $6.9$.
Suppose we wish to achieve $D = 30$.  Then
Equation~\ref{eqn:pmin:roundoff} requires $p \geq 99.7 + 23.3 + 2.9 =
125.9$.  Using MPF-128 will suffice.

The fifth column of Table~\ref{tab:precision} shows lower bounds on
the decimal precision for weighted model counting, according to Equation~\ref{eqn:precision:wmc},
assuming $n = 10^7$ and using rescaling.  We can see that even the precision provided by IEEE Double and ERD
guarantees decimal precisions $8.11$.
Using MPF-128 guarantees decimal precision $30.69$.
We can guarantee these levels of precision 
even when 
performing billions of operations to
compute the weighted model count of a formula with 10 million variables.
Importantly, these bounds hold regardless of the weight assignment, as long as all weights are nonnegative.

\subsection{Generalizing to Other Representations}
\label{sect:ddnnf}

The bound of Equation~\ref{eqn:precision:wmc} applies specifically to decision-DNNF
formulas.  Having a decision variable associated with each disjunction
bounds the depth of the sum operations in a formula to $n$.  Many
decision diagrams with a binary branching structure, including Free
Binary Decision Diagrams (FBDDs)~\cite{wegener:siam:2000}
(a generalization of Ordered BDDs~\cite{bryant:ieeetc:1986,knuth:bdd:2011}), and Zero-suppressed Decision
Diagrams (ZDDs)~\cite{minato:sttt:2001,minato:ijcai:2007} can be
translated into smooth decision-DNNF formulas with a size expansion of
most $n$, and hence the bound of Equation~\ref{eqn:precision:wmc}
holds for these.

For more general d-DNNF formulas, and for decision diagrams with
nonbinary branching structures, including multi-valued decision
diagrams (MDDs)~\cite{srinivasan:iccad:1990} and sentential decision
diagrams (SDDs)~\cite{darwiche:ijcai:2011}, it is difficult to find a
useful error bound that applies to entire classes of formulas.
Instead, given a formula to evaluate,
we propose computing an integer-valued error bound based on the
structure of the formula, and using this bound to guide the selection of
the precision $p$ used in a floating-point evaluation.

We can see with all of these representations that the core requirement
is to compute a rational value $V$ by evaluating an arithmetic
expression $\psi$ consisting of rational constants, products, and
sums, where some of the product and sum operations may have more than
two arguments.  Computing this value with floating-point operations
having precision $p$ will yield an approximation $\approxV$ to the
true value $V$.  We can recursively compute an integer bound $e(\psi)$, such that
$\aerror[\approxV, V] \leq e(V) \,\roundepsilon$.  We assume we
can convert each constant $v$ into its floating-point representation
$\approxv$ with at most one rounding, and therefore $e(v) = 1$.  For a
product of the form $\psi' = \prod_{1 \leq i \leq k} \psi_i$, we can
recursively compute $e(\psi') = 2(k-1) + \sum_{1 \leq i \leq k}
e(\psi_i)$.  Here, the multiplications can be performed via any
association.
For a sum of the form
$\psi' = \sum_{1 \leq i \leq k} \psi_i$, we can compute
$e(\psi') = \lceil \log_2(k-1)\rceil + \max_{1\leq i \leq k} e(\psi_i)$.  Here, the sums should be performed as a balanced tree of binary additions.

For expression $\psi$, we can use the computed bound $e(\psi)$ to
select a fraction size $p$ that guarantees a desired level of
precision.  That is, we require $p \geq 2\,\log_2 e(\psi)$, and to
achieve decimal precision $D$, we require $p \geq D \cdot \log_2 10 +\log_2 e(\psi)$.

\subsection{Experimental Validation}

\begin{figure}
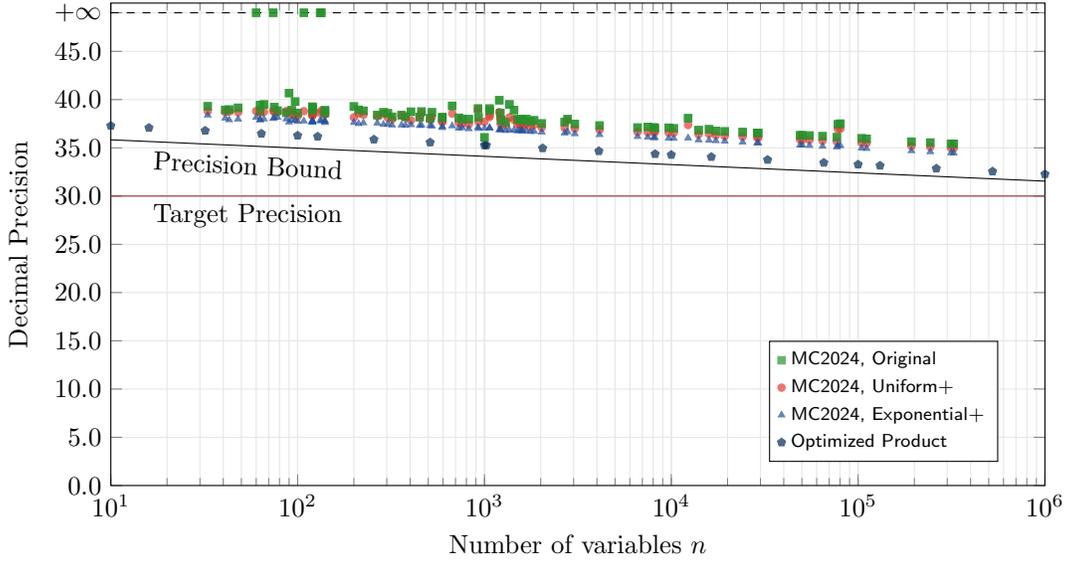

\centering{%
\begin{tikzpicture}
  \begin{axis}[mark options={scale=1.0},height=8cm,width=14cm,grid=both, grid style={black!10}, 
      legend style={at={(0.95,0.30)}},
      legend cell align={left},
                              xmode=log,xmin=10,xmax=1e6,
                              xtick={1,10,100,1000,1e4,1e5,1e6,1e7}, xticklabels={1, $10^1$, $10^2$, $10^3$, $10^4$, $10^5$, $10^6$, $10^7$},
                              ymode=normal,ymin=0, ymax=50,
                              ytick={0, 5, 10, 15, 20, 25, 30, 35, 40, 45, 49},
                              yticklabels={0.0, 5.0, 10.0, 15.0, 20.0, 25.0, 30.0, 35.0, 40.0, 45.0, $+\infty$},
                              xlabel={Number of variables $n$}, ylabel={Decimal Precision}
            ]

    \input{data-formatted/original-mpf+vars}
    \input{data-formatted/upos-mpf+vars}
    \input{data-formatted/epos-mpf+vars}
    \input{data-formatted/optimized-product}
    \legend{
      \scriptsize \textsf{MC2024, Original},
      \scriptsize \textsf{MC2024, Uniform$+$},
      \scriptsize \textsf{MC2024, Exponential$+$},
      \scriptsize \textsf{Optimized Product}
    }
    \addplot[mark=none, dashed] coordinates{(10,49) (1e6, 49)};) 
    \input{data-formatted/original-mpf+vars}
    \addplot[mark=none] coordinates{(1,36.69) (1e7,30.69)};
    \node[right] at (axis cs: 15, 33.2) {\rotatebox{-2.8}{Precision Bound}};
    \addplot[mark=none, color=darkred] coordinates{(1,30) (1e7,30)};
    \node[right] at (axis cs: 15, 28.0) {Target Precision};
 \end{axis}
\end{tikzpicture}
} 
\caption{Decimal Precision Achieved by MPF-128 for Benchmarks with Nonnegative Weights.  The precision is guaranteed to be greater than the bound.
We set as a target to have decimal precisions of at least 30.0.}
\label{fig:pos:mpf}
\end{figure}

To experimentally test the bound of Equation~\ref{eqn:precision:wmc},
we evaluated 200 benchmark formulas from the public and private portions
of the 2024 Weighted Model Counting Competition.\footnote{\url{https://mccompetition.org/2024/mc_description.html}}
We ran version 2 of the
D4 knowledge compiler\footnote{Available at \url{https://github.com/crillab/d4v2}}
to convert these into decision-DNNF\@.
We were able to compile 100 of them within a time
limit of 3600 seconds per formula on a machine with 64~GB of random-access memory.  D4 required a total of 3.82 hours to compile the 100 formulas.

Define a problem \emph{instance} to be a combination of a formula plus
a weight assignment for all of its literals, and a \emph{collection}
as a set of instances, containing multiple formulas, with one or more
weight assignment per formula.  As one collection, we computed the
weighted model count for each compiled formula based using the weight
assignment from the competition.  We refer to this as the
\textsf{Original} collection.  We also generated two collections,
consisting of the compiled formulas with five randomly generated
weight assignments for each formula:
\begin{itemize}
\item \textsf{Uniform$+$}: For each variable $x$, weight $w(x)$ is represented by a 9-digit decimal number selected uniformly in the range
  $\interval{10^{-9},\,1-10^{-9}}$. The weight for $\obar{x}$ is then set to
  $w(\obar{x}) = 1-w(x)$.  Such a weight assignment is typical of those used in recent weighted model counting competitions~\cite{fichte:jea:2020}.
\item \textsf{Exponential$+$}: For each variable $x$, weights $w(x)$ and $w(\obar{x})$
  are drawn independently from an exponential distribution in the range
  $\interval{10^{-9},\,10^{+9}}$.  Each weight is represented by a decimal number with 9 digits to the right of the decimal point.
\end{itemize}


For each instance, we evaluated
the weighted count using MPF-128
to get an approximate weight $\approxw$ and
using MPQ to get an exact weight $w$.  We then evaluated the decimal precision according to Equation~\ref{eqn:digitprecision}.   Our implementation used rescaling for all variables
with $w(x) + w(\obar{x}) \not \in \{0, 1\}$.
Although not required for the instances used in this evaluation,
it will insert smoothing terms when $w(x) + w(\obar{x}) = 0$ for variable $x$.

In addition, we evaluated formulas of the form $\bigwedge_{1\leq i
  \leq n} x_i$ for values of $n$ ranging up to one million using a
single weight for every variable.  For each value of $n$, we swept a parameter space of
weights of the form $1 + 10^{-9}\,k$ for $1 \leq k \leq 1000$ and
chose the value of $k$ that minimized the decimal precision.  We refer
to this as the \textsf{Optimized Product} collection.

Figure~\ref{fig:pos:mpf} shows the result of these evaluations for
the four collections.  For the two collections with multiple weight assignments per formula,
we show only the minimum precision achieved for each formula.
Each data point
represents one combination of formula and weight selection method and is placed
along the X~axis according to the number of variables
and on the Y~axis according to the computed decimal precision.
The plot also shows the precision bound of Equation~\ref{eqn:precision:wmc} for $c=\log_{10} 7$.
Results are shown for 98 of the 100 formulas, since the evaluation consistently ran out of memory when using rational arithmetic for two of the formulas.

As  expected, all data points stay above the precision bound.
Indeed, most exceed the bound by several decimal digits.  Our bound
assumes that rounding either consistently decreases or consistently
increases each computed result.  In practice, rounding goes in both
directions, and therefore the computed results stay closer to the true
values.  The optimized products demonstrate that particular
combinations of formula and weight assignment can come within one
decimal digit of the precision bound and also to track its general
trend.  Indeed, we can use $c = \log_{10} 3$ for these formulas,
since the weighted model count is the product of literal weights.
For $n=10^6$ we get a bound of $32.055$.
Using $w=1.000000453$, we get
a computed decimal precision of $32.273$, a difference from the bound of just $0.218$.

We can see that the achieved decimal precision for the original weight
assignment is somewhat better than for the collections with five
instances per formula.  This can be attributed, in part, to the fact
that plotted values for the other collections show the minimum digit
precision for five weight assignments.  There are even five instances
with the original weight assignments where the values computed with
floating-point arithmetic are exact.  A deeper examination shows these
are particularly simple instances for weighted model counting, with
only 2--3 variables having nonunit weights, and with the counts for
the formulas depending only on the property that the weight for each
of these variables and its complement sum to one.  Importantly, the
data for all collections shows the general trend of the digit
precision decreasing linearly along the logarithmically-scaled X~axis.

As shown in Figure~\ref{fig:pos:mpf}, we select a target precision of
$D=30$ when using floating-point representations with $p \geq 128$.
This target is achieved for our benchmarks with nonnegative weights,
and it should suffice for most applications.

\subsection{Achieving Different Target Precisions}
\label{sect:precision:nonnegative}

\begin{figure}
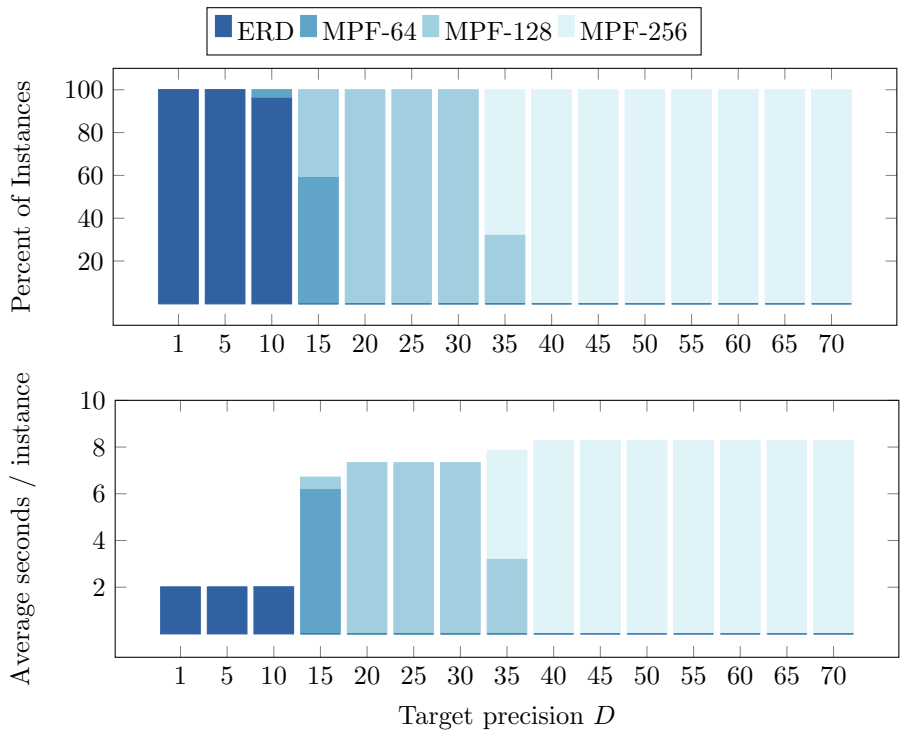

\begin{center}
\begin{tikzpicture}
  \begin{axis}[
      ybar stacked,
      width = 12cm,
      height=5cm,
      bar width=15pt,
      legend style={at={(0.435,1.05)},
        anchor=south,
        legend columns=-1
      },
      ylabel = {Percent of Instances},
      ytick = {20, 40, 60, 80, 100},
      yticklabels ={20, 40, 60, 80, 100},
      symbolic x coords={1, 5, 10, 15, 20, 25, 30, 35, 40, 45, 50, 55, 60, 65, 70},
      xtick=data,
      ]
\input{data-formatted/tabulate-percent-nonneg}
\legend{\strut ERD, \strut MPF-64, \strut MPF-128, \strut MPF-256}
  \end{axis}
\end{tikzpicture}

\begin{tikzpicture}
  \begin{axis}[
      ybar stacked,
      width = 12cm,
      height= 5cm,
      bar width=15pt,
      ylabel = {Average seconds / instance},
      xlabel = {Target precision $D$},
      ymax = 10,
      ytick = {2,4,6,8,10},
      yticklabels = {2,4,6,8,10},
      symbolic x coords={1, 5, 10, 15, 20, 25, 30, 35, 40, 45, 50, 55, 60, 65, 70},
      xtick=data,
      ]
\input{data-formatted/tabulate-aeffort-nonneg}
  \end{axis}
\end{tikzpicture}
\end{center}
\caption{Percent of instances (top) and average evaluation time (bottom) to achieve target precision $D$,
  according to the floating-point representation used.
The evaluation was performed for 1000 instances, all with nonnegative weights.}
\label{fig:solution:nonnegative}
\end{figure}

Figure~\ref{fig:solution:nonnegative} summarizes the performance of
our evaluation strategy for the 1000 instances with nonnegative weights to achieve target
precisions $D$ ranging from 1 to 70.  For these, the minimum fraction
size $p$ is selected according to Equations~\ref{eqn:pmin:linear} and
\ref{eqn:pmin:roundoff}, and the formulas are evaluated using the
floating-point representation providing that level of precision.
These equations depend on the number of variables in the formula, and so some
instances can use lower precision representations than is implied by
the fifth column of Table~\ref{tab:precision}.

The upper part of the figure shows which representations are used for each target precision $D$.
For $D=1$ and $D=5$, all can be evaluated with ERD, as can most of the instances for $D=10$.
MPF-64 suffices for the remaining instances with $D=10$ and the majority with $D=15$.  MPF-128 then suffices through $D=30$,
but achieving higher values of $D$ requires using MPF-256 for most ($D=35$) and then all cases.

\begin{table}
  \caption{Average time to evaluate nonnnegative instances for different numeric representations.
  Double-precision evaluation could fail due to underflow or overflow.  MPQ could fail due to memory limitations.}
  \label{tab:nonneg:perf}
\begin{center}
    \begin{tabular}{llrrrrrr}
      \toprule
Item &        Double & ERD & MPF-64 & MPF-128 & MPF-256 & MPQ \\
      \midrule
Average Seconds & 1.87 & 2.01 & 6.66 & 7.33 & 8.28 & 182.23 \\
Relative to ERD & 0.93$\times$ & 1.00$\times$ & 3.31$\times$ & 3.65$\times$ & 4.12$\times$ & 90.66$\times$ \\
      \bottomrule
    \end{tabular}
\end{center}
\end{table}

The lower portion of Figure~\ref{fig:solution:nonnegative} shows the
average evaluation time per instance (in seconds) as a function of
target precision $D$.  Several trends can be seen here, which are
further highlighted in Table~\ref{tab:nonneg:perf}.  This table shows
the average time per instance for the evaluations using the five
different representations.  The evaluation using double-precision
failed for 628 of the evaluations due to underflow and overflow, while
the evaluations using MPQ failed for 50 due to running out of memory.
The times for the failing cases are included in the averages.  We can
see that evaluation using ERD required only $1.07\times$ longer than
with Double while also successfully evaluating all 1000 instances.
Relative to ERD, the times for evaluation using MPF were
$3.3$--$4.1\times$ longer, with a suprisingly low increase with the
precision.  Finally, evaluation using MPQ requires substantially more
time, in part due to the 50 failing cases.  Even considering only the
successful evaluations gives an average of 106.93 seconds per instance, $52.2\times$ longer than for ERD\@.

\section{Mixed Negative and Positive  Weights}
\label{sect:neg}

The analysis of Section~\ref{sect:nonneg} no longer holds when
some literals have negative weights, while others have positive weights.  With a floating-point
representation, summing combinations of negative and positive
values can cause \emph{cancellation}, where arbitrary levels of
precision are lost~\cite{knuth:fp:1981}.  Consider, for example, the computation
$s + T - T$, where $s$ and $T$ are positive floating-point values, with $s \ll T$.  Using
bounded-precision arithmetic, evaluating the sum as $s + (T - T)$ will yield $s$.
Evaluating it as $(s + T) - T$, however, can yield $0$ or some other value that bears little relation to $s$.
Cancellation can also occur when evaluating a sum $s + T - T'$, where $T \approx T'$.

\subsection{Challenging Formulas and Weight Assignments}

Cancellation can arise when evaluating decision-DNNF formulas to such a degree that no floating-point precision $p$ that grows sublinearly with $n$ will suffice.
As an example, consider the following smooth, decision-DNNF formula $\tau_n$ over $n+1$ variables:
\begin{eqnarray}
\tau_n  & = & z \land \left[\bigwedge_{i = 1}^{n} \obar{x}_i \; \lor \; \bigwedge_{i = 1}^{n} x_i\right] \quad \lor \quad \obar{z} \land \left [\bigwedge_{i = 1}^{n} x_i\right] \label{eqn:max:precision}
\end{eqnarray}
with a weight assignment having only a single literal assigned a negative weight:
\begin{displaymath}
\begin{array}{lllll}
\makebox[1cm]{} &  w(z) \; = \; +1 & \makebox[2cm]{} &  w(x_i) \; = \; 10^{+9} & 1 \leq i \leq n \\
\makebox[1cm]{} &  w(\obar{z}) \; = \; -1 & &  w(\obar{x}_i) \; = \; 10^{-9} & 1 \leq i \leq n \\
\end{array}
\end{displaymath}
Computing $w(\tau_n)$  evaluates the sum $(s + T) - T$, where
$s = 10^{-9n}$ and $T = 10^{+9n}$.  Avoiding cancellation requires using a floating-point representation with a fraction of at least
$p = (18 \cdot \log_2 10)\, n \approx 60 \, n$ bits.
Using MPQ, we were able to compute $w(\tau_{10^7})$ exactly in around 35 seconds, even though the final step requires a total of 1.87~gigabytes to store the arguments
$s+T$ and $-T$, and the result $s$.  In general, however, rational arithmetic can be very time and memory intensive.

\begin{figure}
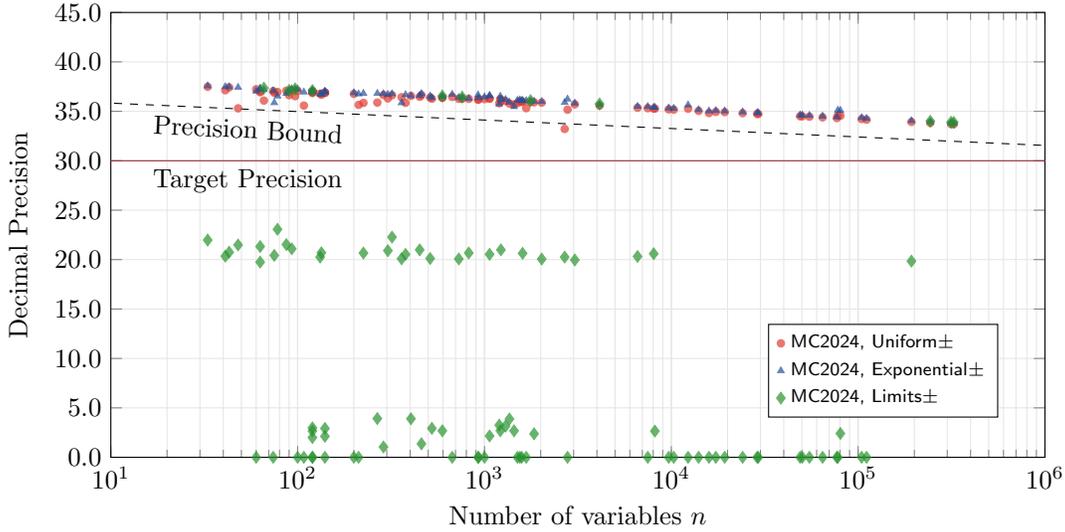

\centering{%
\begin{tikzpicture}
  \begin{axis}[mark options={scale=1.0},height=7.5cm,width=14cm,grid=both, grid style={black!10},
      legend style={at={(0.95,0.30)}},
      legend cell align={left},
                              xmode=log,xmin=10,xmax=1e6,
                              xtick={1,10,100,1000,1e4,1e5,1e6,1e7}, xticklabels={1, $10^1$, $10^2$, $10^3$, $10^4$, $10^5$, $10^6$, $10^7$},
                              ymode=normal,ymin=0, ymax=45,
                              ytick={0, 5, 10, 15, 20, 25, 30, 35, 40, 45, 50},
                              yticklabels={0.0, 5.0, 10.0, 15.0, 20.0, 25.0, 30.0, 35.0, 40.0, 45.0, 50.0},
                              xlabel={Number of variables $n$}, ylabel={Decimal Precision}
            ]

    \input{data-formatted/uposneg-mpf+vars}
    \input{data-formatted/eposneg-mpf+vars}
    \input{data-formatted/bposneg-mpf+vars}
    \legend{
      \scriptsize \textsf{MC2024, Uniform$\pm$},
      \scriptsize \textsf{MC2024, Exponential$\pm$},
      \scriptsize \textsf{MC2024, Limits$\pm$}
    }
    \addplot[mark=none, dashed] coordinates{(1,36.69) (1e7,30.69)};
    \node[right] at (axis cs: 15, 33.2) {\rotatebox{-2.8}{Precision Bound}};
    \addplot[mark=none, color=darkred] coordinates{(1,30) (1e7,30)};
    \node[right] at (axis cs: 15, 28.0) {Target Precision};
 \end{axis}
\end{tikzpicture}
} 
\caption{Decimal Precision Achieved by MPF-128 for Benchmarks with Mixed Weights.  There is no guaranteed bound for precision, but many cases remain above the nonnegative weight bound.
The Limits$\pm$ weight assignment is designed to maximize precision loss.}
\label{fig:posneg:mpf}
\end{figure}

Contrary to the example of
Equation~\ref{eqn:max:precision},
floating-point arithmetic performs surprisingly well for
many real-world
problems, even in the presence of negative weights.
In Figure~\ref{fig:posneg:mpf}, we see a similar plot to that of Figure~\ref{fig:pos:mpf} for
three collections of weight assignments with mixed negative and positive weights.  The first two are generalizations of those used earlier:
\begin{itemize}
\item \textsf{Uniform$\pm$}: For each variable $x$, weights $w(x)$ and $w(\obar{x})$ have magnitudes drawn independently
from a uniform distribution in the range  $\interval{10^{-9},\,1-10^{-9}}$ and are represented as 9-digit decimal numbers.  Each is negated with probability $0.5$.
\item \textsf{Exponential$\pm$}: For each variable $x$, weights $w(x)$ and $w(\obar{x})$ have magnitudes
  drawn independently from an exponential distribution in the range $\interval{10^{-9},\,10^{+9}}$ and are represented with 9 digits to the right of the decimal point.  Each is negated with probability $0.5$.
\end{itemize}
As can be see with these plots, the results mostly stay above the precision bound of Equation~\ref{eqn:digitprecision},
even though this bound need not hold.  All stay above the target precision of $30.0$.

On deeper inspection, we can see that setting up the conditions for a
cancellation of the form $(s + T) - T'$, where $T \approx T'$, requires
1) a large dynamic range among the computed values to give widely different values $s$ and $T$, and 2) sufficient
homogeneity in the computed values that we get two values $T$ and
$T'$ such that $T \approx T'$.  A uniform distribution has neither of
these properties.  An exponential distribution has a large dynamic
range, but the computed values tend to be very heterogenous.

To increase the likelihood of precision loss due to cancellation, we devised the following strategy for generating weight assignments:
\begin{itemize}
\item\textsf{Limits$\pm$}:  For variable $x$, each weight $w(x)$ and $w(\obar{x})$ has a magnitude, chosen at random, of either $10^{-9}$ or $10^{+9}$, and it is set negative with probability $0.5$.
  However, we exclude assignments with $w(x) + w(\obar{x}) = 0$.
\end{itemize}
The idea here is to give large dynamic ranges plus a high degree of
homogeneity.  The plots for this assignment in
Figure~\ref{fig:posneg:mpf} demonstrate the success of this strategy,
with many results falling below the target precision of $30.0$.
This figure
 presents a pessimistic
perspective for these instances, since it only shows the minimum
precision achieved out of five instances in a collection for each formula.  Considering all 490
instances,
221 (45\%) yielded results above the
target precision of $30.0$.  We can also see how our choice of weights
leads to two bands of low precision.  129 instances (26\%) had digit
precisions in a band between $19.0$ and $23.3$.  These were ones where
the evaluation encountered values of $s$ and $T$ that
differ by a factor of around $10^{18}$.  The remaining 140 instances
(29\%) had decimal precisions below $5.0$.  These were ones where the
encountered values of $s$ and $T$ differed by a factor of around
 $10^{36}$.

\subsection{Interval Computation Applied to Weighted Model Counting}
\label{sect:interval}

\begin{figure}
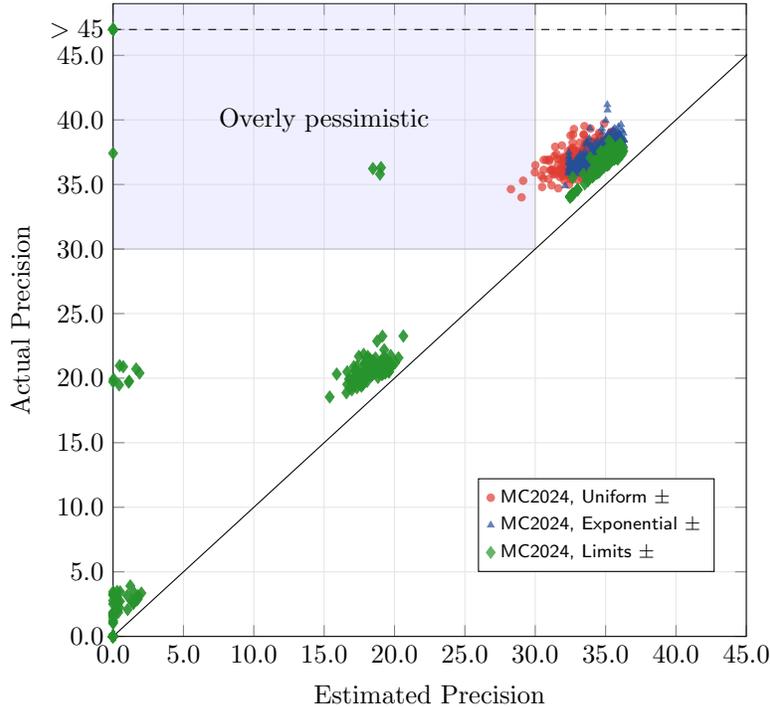

\centering{%
\begin{tikzpicture}
  \begin{axis}[mark options={scale=1.0},height=10cm,width=10cm,grid=both, grid style={black!10},
      legend style={at={(0.95,0.25)}},
      legend cell align={left},
                              xmode=normal,xmin=0, xmax=45,
                              xtick={0, 5, 10, 15, 20, 25, 30, 35, 40, 45, 50},
                              xticklabels={0.0, 5.0, 10.0, 15.0, 20.0, 25.0, 30.0, 35.0, 40.0, 45.0, 50.0},
                              ymode=normal, ymin=0, ymax=49,
                              ytick={0, 5, 10, 15, 20, 25, 30, 35, 40, 45, 47},
                              yticklabels={0.0, 5.0, 10.0, 15.0, 20.0, 25.0, 30.0, 35.0, 40.0, 45.0, $> 45$, 50.0},
                              xlabel={Estimated Precision}, ylabel={Actual Precision}
            ]

    \draw[ fill={blue!30}, opacity=0.2] (axis cs: 0,30) rectangle (axis cs: 30,50);
    \node at (axis cs:15,40) {Overly pessimistic};

    \input{data-formatted/uposneg-mpfi-est+act}
    \input{data-formatted/eposneg-mpfi-est+act}
    \input{data-formatted/bposneg-mpfi-est+act}
    \legend{
      \scriptsize \textsf{MC2024, Uniform $\pm$},
      \scriptsize \textsf{MC2024, Exponential $\pm$},
      \scriptsize \textsf{MC2024, Limits $\pm$}
    }
    \addplot[mark=none, dashed] coordinates{(0,47) (45,47)};
    \addplot[mark=none] coordinates{(0,0) (50,50)};
    \input{data-formatted/bposneg-mpfi-est+act}
 \end{axis}
\end{tikzpicture}
} 
\caption{Predictive Accuracy of MPFI-128 Interval Arithmetic.  MPFI never has a higher estimate than the actual, but it can
incorrectly predict a precision less than the target of 30.}
\label{fig:mpfi}
\end{figure}

We can see from Figure~\ref{fig:posneg:mpf} that floating-point
evaluations generate accurate results in many cases, but we must be
able to discern when those occur.  Given that capability, we can
devise a strategy that combines multiple methods to reliably compute
weighted counts.  We use a target precision bound of $D=30$ here for
illustrative purposes.

Interval arithmetic provides a mechanism for using the approximate
computations of floating-point arithmetic, while providing a
guaranteed precision for the result.  It will only be beneficial,
however, if the interval bounds remain tight enough that the digit
precision bound of Equation~\ref{eqn:interval:digitprecision} meets
our target decimal precision.  Our target bound of 30 seems fairly
aggressive in this respect: the width of the interval $\vmax-\vmin$
must be over 30 orders of magnitude smaller than the magnitudes of
$\vmin$ and $\vmax$.  Even the instances with only nonnegative weights
had decimal precisions as low as $34.5$, and so there is not much
room for further degradation.

Figure~\ref{fig:mpfi} shows the result of evaluating 100 formulas for
the three weight assignment collections containing mixed weights,
with five instances per collection for each formula.  
The evaluation uses
MPFI, with $p=128$ (we refer to this as ``MPFI-128'') to get an estimated decimal precision (X~axis) and
a nominal weight (the midpoint of the interval), along with MPQ to get
the exact weight.  The actual precision (Y~axis) is computed based on the
nominal and actual weights.  The evaluations using MPQ consistently runs out of
memory for two of the formulas, and hence the plot shows 1470
data points.  
Every point lies above the diagonal line where the two precisions are equal---the interval computation never overestimates the digit
precision.

Overall, we can see that the interval estimates are quite reliable,
especially for predicting which computed weights exceed the target
threshold of 30.
The
interval computations determines that 1189 ($80.9\%$) instances are
above the target threshold: 490 from \textsf{Exponential$\pm$}, 486 from
\textsf{Uniform$\pm$}, and 213 from \textsf{Limits$\pm$}.
Points lying in the blue rectangle
indicate instances where the estimate is overly pessimistic: they
estimate a target precision below 30, while the actual precision is
above.  This occurs for only 20 of the 1470 instances ($1.4\%$).  Of
these, 4 are from the \textsf{Uniform$\pm$} collection, while 16 are
from the \textsf{Limits$\pm$} collection.  

The interval analysis captures the general trend shown in
Figure~\ref{fig:posneg:mpf} that, even with mixed
weights, floating-point evaluation only degrades significantly due to cancellation
for the weight assignments designed to maximize
this effect.  This gives us hope that we can use interval computation
to handle a large portion of instances having mixed weights.

\subsection{Achieving Different Target Precisions}
\label{sect:precision:mixed}

\begin{figure}
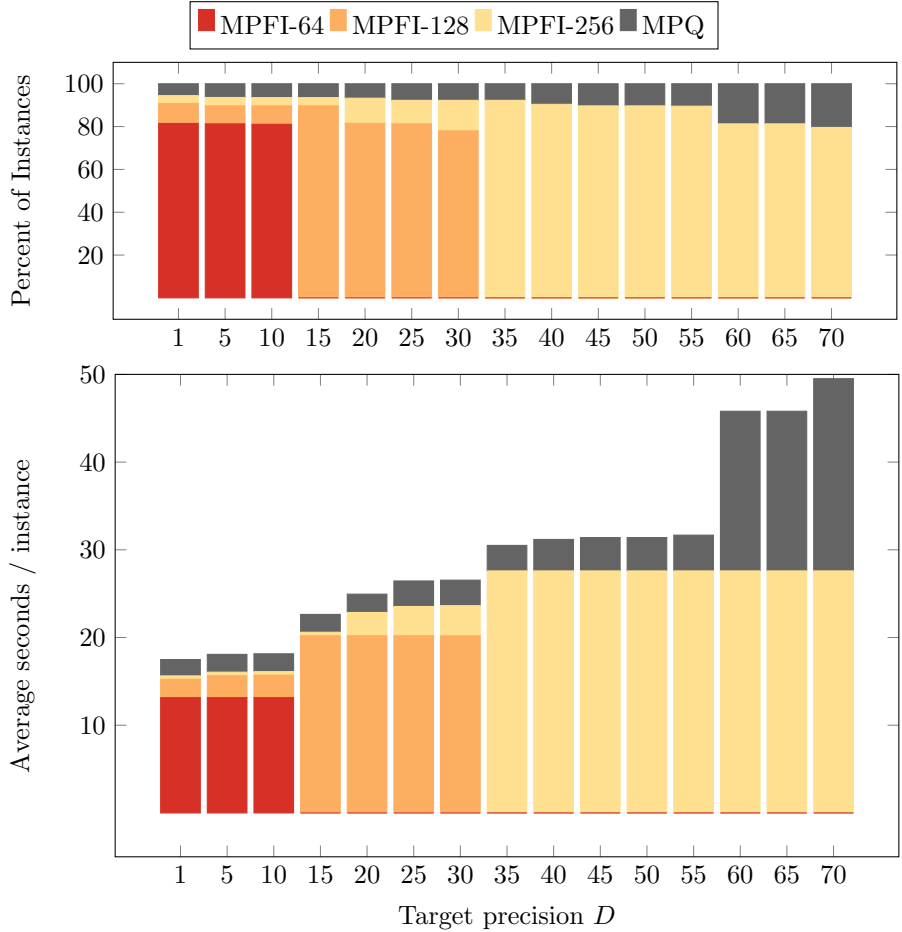

\begin{center}
\begin{tikzpicture}
  \begin{axis}[
      ybar stacked,
      width = 12cm,
      height=5cm,
      bar width=15pt,
      legend style={at={(0.435,1.05)},
        anchor=south,
        legend columns=-1
      },
      ylabel = {Percent of Instances},
      ytick = {20, 40, 60, 80, 100},
      yticklabels ={20, 40, 60, 80, 100},
      symbolic x coords={1, 5, 10, 15, 20, 25, 30, 35, 40, 45, 50, 55, 60, 65, 70},
      xtick=data,
      ]
\input{data-formatted/tabulate-percent-posneg}
\legend{\strut MPFI-64, \strut MPFI-128, \strut MPFI-256, \strut MPQ}
  \end{axis}
\end{tikzpicture}

\begin{tikzpicture}
  \begin{axis}[
      ybar stacked,
      width = 12cm,
      height= 8cm,
      bar width=15pt,
      ylabel = {Average seconds / instance},
      xlabel = {Target precision $D$},
      ymax = 50,
      ytick = {10, 20, 30, 40, 50},
      yticklabels = {10, 20, 30, 40, 50},
      symbolic x coords={1, 5, 10, 15, 20, 25, 30, 35, 40, 45, 50, 55, 60, 65, 70},
      xtick=data,
      ]
\input{data-formatted/tabulate-aeffort-posneg}
  \end{axis}
\end{tikzpicture}
\end{center}
\caption{Percent of instances (top) and average evaluation time (bottom) to achieve target precision $D$,
  according to the numerical representation used.
The evaluation was performed for 1500 instances, with mixed negative and positive weights.}
\label{fig:solution:mixed}
\end{figure}

Figure~\ref{fig:solution:mixed} illustrates the performance of a
simple method for achieving target precisions $D$ ranging from $1$ to
$70$ for the 1500 instances with mixed weight assignments.  For each
formula and target precision, it selects a starting precision based on
Equations~\ref{eqn:pmin:linear} and \ref{eqn:pmin:roundoff}, even
though these bounds are not guaranteed.  It then iterates using MPFI
with increasing levels of precision (64, 128, 256) until the target
precision can be guaranteed.  If all of these fails, it performs the
evaluation with rational arithmetic using MPQ\@.  For example, to
achieve target precision $D=1$, 1222 instances completed with MPFI-64,
140 with MPFI-128, and 55 with MPFI-256.  That left 83 instances to
evaluate using MPQ\@.  Achieving higher target precisions follows the same pattern, such that 306 of the instances require evaluation with MPQ
to reach the target precision of $D=70$.

The lower part of Figure~\ref{fig:solution:mixed} shows the average
time per instance with this approach.  We can see that these times are
significantly larger than those for nonnegative weights (Figure~\ref{fig:solution:nonnegative}), especially
since we have no counterpart to ERD for mixed weights.  We see also
that the evaluations with MPFI tend to have a greater sensitivity to
precision, and that evaluations using MPQ incur a significant performance penalty.

This iterative approach incurs wasted effort when an evaluation using
MPFI fails to achieve the target precision.  Overall,
however, the wasted effort is below $12\%$ of the total execution
time for each of the target precisions.

\section{A Hybrid Approach}
\label{sect:hybrid}

\begin{table}
  \caption{Performance Comparison of Different Implementation
    Strategies for Target Precision $D=30$.  Run entries of the form
    $S$+$F$ indicate that $S$ runs were successful and $F$ runs either
    ran out of memory or failed to meet the target precision.  Our
    hybrid strategy is shown in red.}
  \label{tab:compare}
  \centering{
  \begin{tabular}{llrrrrr}
    \toprule
    \multicolumn{1}{c}{Strategy} & & \multicolumn{1}{c}{MPF-128} & \multicolumn{1}{c}{MPFI-128} & \multicolumn{1}{c}{MPFI-256} & \multicolumn{1}{c}{MPQ} & \multicolumn{1}{c}{Combined}
    \\
    \midrule
    MPQ only & Runs &  &  &  & 2450+50 & 2450+50 \\
 & Hours &  &  &  & 95.22 & 95.22 \\
\midrule
MPF|MPQ & Runs & 1000+0 &  &  & 1470+30 & 2470+30 \\
 & Hours & 2.04 &  &  & 57.15 & 59.19 \\
\midrule
MPF|MPFI$\times$1+MPQ & Runs & 1000+0 & 1215+285 &  & 281+4 & 2496+4 \\
 & Hours & 2.04 & 8.43 &  & 7.80 & 18.26 \\
\midrule
\textcolor{red}{MPF|MPFI$\times$2+MPQ} & Runs & 1000+0 & 1215+285 & 169+116 & 116+0 & 2500+0 \\
 & Hours & 2.04 & 8.43 & 1.43 & 1.21 & 13.10 \\
\midrule
MPF|MPFI-256+MPQ & Runs & 1000+0 &  & 1384+116 & 116+0 & 2500+0 \\
 & Hours & 2.04 &  & 11.50 & 1.21 & 14.75 \\

    \\[-1em]
   \bottomrule
  \end{tabular}
  } 
\end{table}

We can combine our three approaches---floating-point
arithmetic, interval computation, and rational arithmetic---into a
single, hybrid approach.  We consider a target of $D=30$, although the same strategy applies for other target precisions.
We can measure 
measure the effectiveness of our scheme based on 2500 instances---100 formulas, each with five collections of five instances, as shown in the fourth entry in Table~\ref{tab:compare}.
\begin{enumerate}
\item For instances where all weights are nonnegative, use MPF-128,
  relying on the bound of Equation~\ref{eqn:precision:wmc} to
  guarantee sufficient precision.  This evaluation succeeded for all 1000 such
  instances, including 20 for which the evaluation with rational arithmetic failed.
\item For instances with mixed weights, attempt evaluations with increasing precision and cost:
\begin{enumerate}
\item 
  Use MPFI-128.  If the estimated precision
  bound meets the target bound, then we are done.  This succeeded for
  1215 of the 1500 instances evaluated, including 26 for which the evaluation with rational arithmetic failed.
\item For instances where the estimated precision does not meet the target, perform a second run with MPFI-256.  This succeeded
  for 169 of the 285 instances evaluated, including 4 for which the evaluation with rational arithmetic failed.
\item When the second attempt at interval computation fails, evaluate with rational arithmetic using MPQ\@.
This
  succeeded for the remaining 116 instances.
\end{enumerate}
\end{enumerate}
Overall this strategy succeeded for all 2500 instances.

Table~\ref{tab:compare} summarizes the performance of five different strategies, with our hybrid strategy as the fourth.
Evaluating all 2500
instances with MPQ completes
2450 of them, requiring a
total of $95.2$ hours, of which over 
$22$ hours is spent on the 50 failed runs.
Combining MPF-128 for the instances with nonnegative weights with MPQ for the rest
completes an
additional 20 instances and drops the total time to $55.2$ hours.
Using one pass with MPFI-128
for the instances with mixed weights and then using MPQ for those
that do not meet the target precision completes all but 4 instances
and drops the total time to $18.3$ hours.
Our proposed hybrid approach completes all 2500 instances in a total of
$13.1$ hours.
Finally, skipping the MPFI-128 evaluation and instead going directly to MPFI-256 avoids
wasted effort, but that
does not
compensate for the time required to perform all 1500 evaluations 
with $p=256$.  

The impact of the time spent in compilation versus in weighted
evaluation depends on the usage model.  With these benchmarks, the 25
instances for each formula can be evaluated after compiling the
formula once.  Thus, the time for compilation plus evaluation ranges
from $16.9$ hours for the hybrid method to $99.0$ hours when only
using rational arithetic, giving the hybrid method a total speedup of
$5.90\times$.  On the other hand, if we require each instance to be
compiled separately, the total time would range from $108.6$ to
$190.7$ hours, giving a speedup of around $1.76\times$.  Some
applications of weighted counting require many different evaluations
of a single formula~\cite{sundermann:tsem:2024}; these would benefit
the most from improvements in the evaluation speed. Importantly, the
hybrid method completes all 2500 instances, whereas rational
arithmetic fails for 100 of them.

\begin{figure}[t]
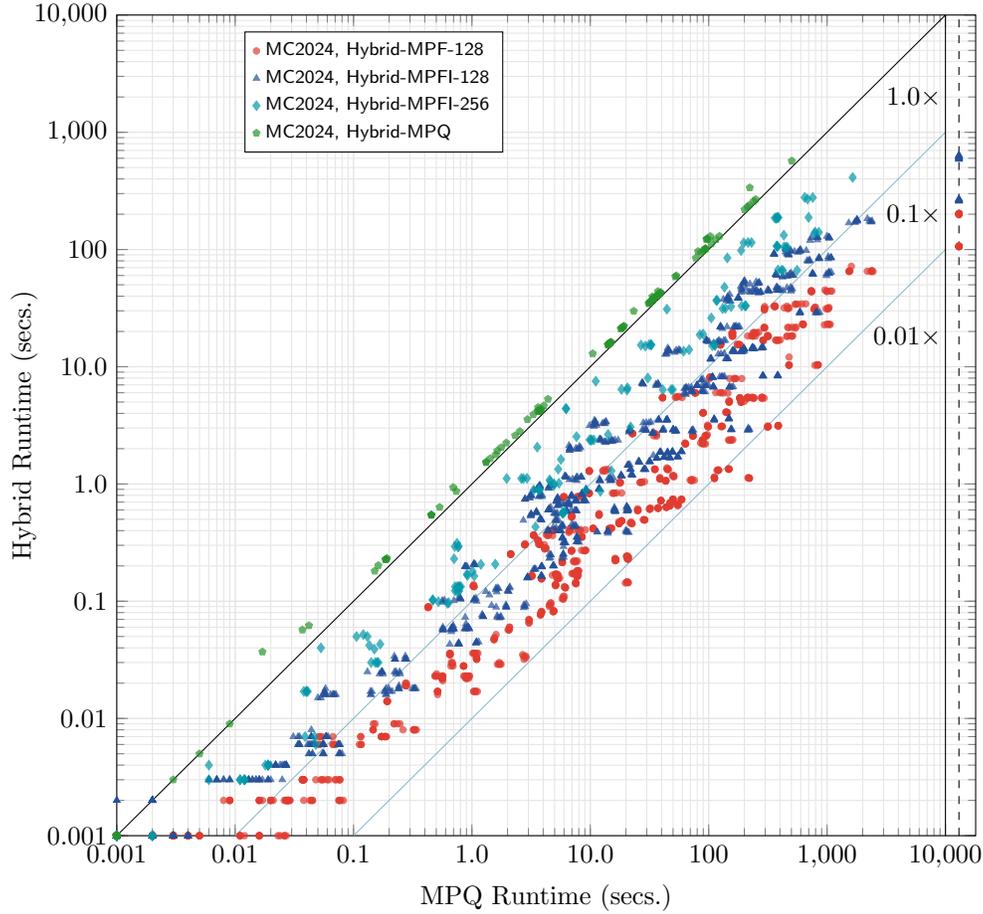

\centering{%
\begin{tikzpicture}
  \begin{axis}[mark options={scale=1.0},height=12.5cm,width=13.0cm,grid=both, grid style={black!10},
      legend style={at={(0.45,0.98)}},
      legend cell align={left},
                              xmode=log,xmin=0.001, xmax=18000.0,
                              xtick={1e-3, 1e-2, 1e-1, 1, 10, 1e2, 1e3, 1e4},
                              xticklabels={$0.001$, $0.01$, $0.1$, $1.0$, $10.0$, $100$, $1{,}000$, $10{,}000$},
                              ymode=log, ymin=0.001, ymax=10000.0,
                              ytick={1e-3, 1e-2, 1e-1, 1, 10, 1e2, 1e3, 1e4},
                              yticklabels={$0.001$, $0.01$, $0.1$, $1.0$, $10.0$, $100$, $1{,}000$, $10{,}000$},
                              xlabel={MPQ Runtime (secs.)}, ylabel={Hybrid Runtime (secs.)}
            ]

    \input{data-formatted/combo-mpf+mpq}
    \input{data-formatted/combo-mpfi+mpq}
    \input{data-formatted/combo-mpfi2+mpq}
    \input{data-formatted/combo-mpq+mpq}
    \legend{
      \scriptsize \textsf{MC2024, Hybrid-MPF-128},
      \scriptsize \textsf{MC2024, Hybrid-MPFI-128},
      \scriptsize \textsf{MC2024, Hybrid-MPFI-256},
      \scriptsize \textsf{MC2024, Hybrid-MPQ}
    }
    \input{data-formatted/combo-mpf+nompq}
    \input{data-formatted/combo-mpfi+nompq}
    \input{data-formatted/combo-mpfi2+nompq}
    \addplot[mark=none] coordinates{(0.001,0.001) (10000.0,10000.0)};
    \node[left] at (axis cs: 11000, 2000) {$1.0\times$};
    \addplot[mark=none, color=lightblue] coordinates{(0.01, 0.001) (10000.0, 1000.0)};
    \node[left] at (axis cs: 11000, 200) {$0.1\times$};
    \addplot[mark=none, color=lightblue] coordinates{(0.1, 0.001)  (10000.0, 100.0)};
    \node[left] at (axis cs: 11000, 18) {$0.01\times$};
    \addplot[mark=none] coordinates{(1e4,1e-3) (1e4,1e4)};
    \addplot[mark=none,dashed] coordinates{(1.3e4,1e-3) (1.3e4,1e4)};
 \end{axis}
\end{tikzpicture}
} 
\caption{Runtime for hybrid method vs.~for MPQ, categorized by the solution method.  
Successful evaluation with MPF or MPFI can significantly reduce the runtime,
but failed evaluations cause some overhead.}
\label{fig:runtime}
\end{figure}

Figure~\ref{fig:runtime} compares the runtimes  for the hybrid
strategy (Y~axis) with target precision $D=30$, versus that for performing an evaluation using rational
arithmetic (X~axis) for all 2500 instances.
These are categorized by the method by which the hybrid
method completed.  The diagonal lines show the relative time for the
hybrid approach versus rational arithmetic.
The points on the right indicate the 50 instances where the evaluation using MPQ fails, but the hybrid method completes.

Of the 2450 instances where the MPQ evaluation completed, the 980 with
nonnegative weights can be evaluated using MPF-128.  Many of these
also have very small runtimes, even for MPQ\@.  Considering just the
670 instances for which MPQ requires more than $1.0$ seconds, we find
that MPF-128 runs between $7.4$ and $197.6$ times faster than MPQ, with an
average of $38.1$ and a median of $32.4$.  This shows a clear
performance benefit in using MPF when all weights are nonnegative.

Of the 1470 instances  containing mixed weights where the MPQ evaluation completes, 1189 ($80.9\%$)
are successfully evaluated using MPFI-128.  Considering the
824 instances for which MPQ requires more than $1.0$ seconds, we find
that MPFI-128 runs between $3.1$ and $75.8$ times faster, with an average of
$14.2$ and a median of $11.4$.  Again, this level of evaluation has a
clear performance benefit.  An additional 165 instances
($11.2\%$) are successfully evaluated using MPFI-256.  Of
the 97 instances for which MPQ requires more than $1.0$ seconds, we
find that the combined time for two runs with MPFI range between $1.4$ and $14.0$ times faster,
with an average of $4.6$ and a median of $4.1$.  This level of evaluation also provides a performance benefit.
Finally, 116 instances ($7.9\%$) require an
evaluation using MPQ\@.  In these cases, the hybrid runtime is
greater than that for MPQ alone, since the program also performs two
evaluations using MPFI\@.  Of the 81 instances for which MPQ requires
more than $1.0$ seconds, we find that the hybrid approach runs between
$1.05$ and $1.51$ times slower, with an average of $1.14$ and a median of $1.12$.
Fortunately, this performance penalty is more than offset by the gains achieved by the less costly evaluation methods.


\section{The Extended-Range Double (ERD) Floating-Point Representation}
\label{sect:erd}

As observed in Section~\ref{sect:background:numbers}, the 11-bit exponent
  field of the IEEE Double representation limits the range of
  representable numbers (excluding infinities) to around
  $10^{\pm308}$~\cite{overton:siam:2001}.  We can overcome this
  limitation by representing floating-point numbers as a pair
  $\langle d, e\rangle$, where $d$ is a floating-point number in IEEE double
  format, and integer exponent $e$ is represented as a 64-bit signed
  integer.
  Adding this exponent field greatly expands the representable range
  of numbers.  We can do so while having the hardware support for
  double-precision arithmetic take care of the trickiest parts of conversion, addition, and multiplication.

  For most IEEE double values $d$, the exponent value $\fexp(d)$
  has a range $-1022 \leq \fexp(d) \leq +1023$.\footnote{The exponent is stored in \emph{biased} form~\cite{overton:siam:2001},
  but for our presentation we considered its unbiased value.}
  The fraction value
  $\fexp(d)$ satisifes $1.0 \leq \ffrac(d) < 2.0$.
  Value $0.0$ is stored with a
  special exponent value, as are denormalized numbers, infinities,
  and not-a-number (NaN).  We do not support the latter three cases
  with ERD\@.

  We will say that the pair $\langle d, e \rangle$ is
  \emph{normalized} when either $d = 0.0$ and $e = 0$, or $d$ is
  nonzero, but it has an exponent value $\fexp(d) = 0$.  An arbitrary
  pair $\langle d, e \rangle$ can be normalized as $\langle 0.0, 0 \rangle$
  when $d = 0$, or as $\langle d', e + \fexp(d)\rangle$ when
  $d \not = 0$, where $d'$ has the same sign and fraction as $d$, but
  an exponent value of $0$.

  Multiplying a set of ERD values of the form
$\langle d_i, e_i \rangle$ for $1 \leq i \leq n$
  can be performed by computing $d = \prod_{1\leq i \leq n} d_i$
  and $e = \sum_{1 \leq i \leq n} e_i$ and then normalizing
  the pair $\langle d, e \rangle$.  Note, however, that $d$
  has the possible range $1.0 \leq d < 2^n$, and so overflow can occur for $n > 1023$.
  Products of longer sequences can be computed by normalizing intermediate results.
  
  Adding a pair of ERD values of the form $\langle d_1, e_1 \rangle$
  and $\langle d_2, e_2 \rangle$ requires considering individual
  cases.  When $d_1 = 0.0$ (respectively, $d_2 = 0.0$),
  the result will be
  $\langle d_2, e_2 \rangle$ (resp., $\langle d_1, e_1\rangle$).
  When $e_1 > 54 + e_2$ (respectively, $e_2 > 54 + e_1$)
  the result will be $\langle d_1, e_1 \rangle$
  (resp., $\langle d_2, e_2 \rangle$).
  Otherwise, we
  normalize the pair $\langle d_1' + d_2, e_2 \rangle$, where $d_1'$
  has the same sign and fraction as $d_1$, but it has exponent
  $e_1 - e_2$.

  We could extract values from and insert values into the exponent
  field of a double-precision number using library functions
  \texttt{frexp} and \texttt{ldexp}~\cite{jones:1991}, but we obtained better
  performance, using our own bit-manipulation code.  The
  compiler was able to optimize the generated machine code with these bit manipulations
  using inline substitution.

\section{Conclusions}
\label{sect:conclusion}

For many applications, floating-point arithmetic can introduce
significant errors due to rounding, and it does not provide any way to
quantify the error.  This paper shows that such uncertainty can be
avoided for weighted model counting.  When all weights are
nonnegative, results can be computed using floating point with
guaranteed precision.  When some weights are negative, the program can
attempt one or more levels of interval computation, and these should
handle a large fraction of the instances.  Ultimately, the program may
need to use rational arithmetic, but the number of such cases should
be small.  By including formulas and weight assignments that are
especially challenging from a numerical perspective in our evaluations, we can be confident
of the robustness of our approach.

\newpage
\bibliography{references}

\end{document}